\documentclass[10pt]{article}

\usepackage{amssymb} \usepackage{amsthm}
\usepackage{latexsym} \usepackage{exscale}
\usepackage{epsfig}
\usepackage{epsf} \usepackage{graphicx}
\newtheorem{theorem}{Theorem}[section]
\newtheorem{lemma}[theorem]{Lemma}

\newtheorem{definition}[theorem]{Definition}

\makeatletter \@addtoreset{equation}{section} \makeatother

\newcommand{\beq}{\begin{equation}}
\newcommand{\eeq}{\end{equation}}

\newcount\gnum
\def\G#1{
  \gnum=`#1
  \ifnum \gnum>`Z
     \advance\gnum by-`a
     \ifcase \gnum \alpha \or \beta \or \gamma \or \delta \or
        \epsilon\or\phi\or\theta\or\eta\or\iota\or\or\kappa\or
        \lambda\or\mu\or\nu\or\or\pi\or\xi\or\rho\or\sigma\or
        \tau\or\upsilon\or\or\omega\or\chi\or\psi\or\zeta\fi
  \else
     \advance\gnum by-`A
     \ifcase \gnum \Alpha \or \Beta \or \Gamma \or \Delta \or
        \Epsilon\or\Phi\or\Theta\or\Eta\or\Iota\or\or\Kappa\or
        \Lambda\or\Mu\or\Nu\or\or\Pi\or\Xi\or\Rho\or\Sigma\or
        \Tau\or\Upsilon\or\or\Omega\or\Chi\or\Psi\or\Zeta\fi
  \fi}

\newtheorem{cor}{Corollary}

\def\({\left(\begin{array}{cccccc}}
\def\){\end{array}\right)}

\def\bes{\begin{eqnarray}}
\def\ees{\end{eqnarray}}

\begin{document}
\title{Formation of singularity and smooth wave propagation
for the non-isentropic compressible Euler equations}

\author{Geng Chen\thanks{chen@math.psu.edu}\thanks{This
paper is part of the first author's Ph.D dissertation at the
University of Massachusetts, Amherst.}\vspace{.2cm}
\\
Department of Mathematics, Penn State University\\
University Park, PA 16802}

\date{}

\maketitle

\begin{abstract}
We define compressive and rarefactive waves and give the
differential equations describing smooth wave steepening for the
compressible Euler equations with a varying entropy profile and
general pressure laws. Using these differential equations, we
directly generalize P. Lax's singularity (shock) formation results
in \cite{lax2} for hyperbolic systems with two variables to the
$3\times3$ compressible Euler equations for a polytropic ideal gas.
Our results are valid globally without restriction on the size of
the variation of initial data.\end{abstract}


\section{Introduction}
In this paper, we consider the initial value problem for the
compressible Euler equations in Lagrangian coordinates in one space
dimension,
\begin{eqnarray}
\tau_t-u_x=0,\label{lagrangian1}\\
u_t+p_x=0,\label{lagrangian2}\\
(\frac{1}{2}u^2+e)_t+(up)_x=0, \label{lagrangian3}\end{eqnarray}
where  $\rho$ is the density, $\tau=\rho^{-1}$ is the specific
volume, $p$ is the pressure, $u$ is the velocity, and $e$ is the
internal energy. For $C^1$ solutions, (\ref{lagrangian3}) can be
equivalently replaced by \beq S_t=0\label{s con},\eeq where $S$ is
the entropy. When the entropy is a constant (isentropic fluid),
(\ref{lagrangian1}) and (\ref{lagrangian2}) become a complete
system, known as the $p$-system \cite{smoller}.

The formation of shock waves from smooth initial data is one of the
central problems of conservation laws. The formation of shock waves
takes place when one or more gradients blows up. In this paper, we
consider the formation of such kind of singularity from smooth
initial data. In \cite{lax2}, {Lax} gave singularity formation
results for hyperbolic conservation laws with two variables,
including the $p$-system as an example. This says that the
singularity forms from smooth initial data if and only if the
initial data includes compressive waves. For the conservation laws
with a coordinate system of Riemann invariants, similar singularity
formation results hold, c.f. \cite{Dafermos}. However, for the
conservation laws with more than two variables and without a
coordinate system of Riemann invariants, including
(\ref{lagrangian1})$\sim$(\ref{lagrangian3}), the singularity
formation results are only for restrictive initial data, see
\cite{Fritz John}, \cite{Liu1}, \cite{Li daqian}. In this paper, we
directly generalize Lax's singularity formation results for the
$p$-system to (\ref{lagrangian1})$\sim$(\ref{lagrangian3}). The
singularity formations for the compressible Euler equations in
multiple space dimensions are considered by \cite{sideris} and
\cite{Rammaha}.

Comparing to the $p$-system, one of the main difficulties for
studying (\ref{lagrangian1}) $\sim$ (\ref{lagrangian3}), even for
the smooth solutions, is the smoothly varying entropy field. Wave
propagation in a smoothly varying entropy field is still not well
understood: even some basic questions, like what should be
considered as the compressive and rarefactive waves have not
previously been answered. Furthermore, the study of smoothly varying
entropy field is also the basis for studying the formation of
singularity, and interactions between shock waves and simple waves.

In an isentropic domain, people have long known how to define the
compressive and rarefactive waves for smooth solutions by the
changes of density, pressure, wave speed or Riemann invariants, c.f.
\cite{courant}, \cite{lax}, \cite{lax2}, \cite{young blake 1},
\cite{young Global interaction}. However, in a smoothly varying
entropy field, there are no Riemann invariants, and the forward and
backward waves can't be geometrically divided. In order to define
rarefactive and compressive waves, we need to find a variable which
is not impacted by the variation of entropy but still can
discriminate rarefactive and compressive waves. Pressure is an
appropriate variable, since it is invariant inside the stationary
entropy field.

We use superscripts $\prime$ and $\backprime$ to denote the
directional derivatives on the forward and backward characteristic
lines, respectively: \beq{\prime}=\partial_t+c
\partial_x,\ {\backprime}=\partial_t-c \partial_x,\label{prime backprime}\eeq
where $c$ is the wave speed, c.f. \cite{lax2}. We define the
Rarefactive ({R}) and Compressive ({C}) waves (or $R/C$ character)
as follows: {\definition\label{def1}{\emph{Consider a solution of
(\ref{lagrangian1})$\sim$(\ref{lagrangian3}), which is smooth in an
open set $U$ in the $(t,x)$-plane and $A$ is a point in $U$. We say
that the solution is forward (backward) rarefactive at $A$, if and
only if $p^{\backprime}<0$ ($p^{\prime}<0$); it is forward
(backward) compressive at $A$, if and only if $p^{\backprime}>0$
($p^{\prime}>0$). }}}\vspace{.2cm}\\
Here, we use the derivative of pressure in the \emph{opposite}
direction to define the \emph{R/C} character, which helps us
discount the disturbance from the waves in the opposite direction.
This definition gives us appropriate physical explanation of
\emph{R/C} characters. This definition is not only for the simple
waves \cite{smoller}, but also for waves in the wave interaction
regions.

Then we restrict our consideration on the polytropic ideal gas
dynamics, where \beq p=Ke^{\frac{S}{c_v}}\tau^{-\gamma},\eeq $K$ and
$c_v$ are positive constants, and $\gamma>1$ is the adiabatic gas
constant. We introduce variables $z$ and $m$: \beq
z=\frac{2\sqrt{K\gamma}}{\gamma-1} \tau^{-\frac{\gamma-1}{2}},\quad
m=e^{\frac{S}{2c_v}}.\label{intr z m}\eeq

In this paper, the smooth solution means that $u,\tau,S$ are $C^2$,
which is equivalent to that $u,z,m$ are $C^2$ and $z$, $m$ are both
positive, by (\ref{intr z m}). If the initial data are smooth, $m$
is smooth and positive until either the $\tau$ or $u$ profile has
some singularities (gradient blowups), since the entropy is
stationary in $C^1$ solutions by (\ref{s con}). Since we only
consider smooth wave propagation and singularity formation from
smooth initial data in this paper, the entropy profile is always
smooth and stationary.

We introduce variables \beq \left.\begin{array}{l}
\alpha=-\frac{p^{\backprime}}{c^2}=u_x+mz_x+\frac{\gamma-1}{\gamma}m_x
z,\end{array}\right.\eeq \beq
\left.\begin{array}{l}\beta=-\frac{p^{\prime}}{c^2}=u_x-mz_x-
\frac{\gamma-1}{\gamma}m_x z.\end{array}\right.\eeq $\alpha$ and
$\beta$ are the generalization of $s_x$ and $r_x$ in a smoothly
varying entropy field, where $s$, $r$ are Riemann invariants in a
constant entropy field. {\theorem\label{remark} The smooth solutions
in (\ref{lagrangian1})$\sim$(\ref{lagrangian3}) satisfy \beq
\left.\begin{array}{l}\alpha^\prime=k_1\{k_2
(3\alpha+\beta)+\alpha\beta-\alpha^2\}, \label{rem1}
\end{array}\right.\eeq
and \beq \left.\begin{array}{l}\beta^\backprime=k_1\{-k_2
(\alpha+3\beta)+\alpha\beta-\beta^2\},\label{rem2}
\end{array}\right.\eeq
where \beq\left.\begin{array}{l}
k_1=\frac{(\gamma+1)K_c}{2(\gamma-1)} z^{\frac{2}{\gamma-1}}, \quad
k_2=\frac{\gamma-1}{\gamma(\gamma+1)}z m_x, \quad K_c\ is\ a\
positive\ constant.\label{k def}\end{array}\right.\eeq}
\\
The equations (\ref{rem1}) and (\ref{rem2}) are not pure ODEs
because they aren't closed.

Then we transform (\ref{rem1}) and (\ref{rem2}) into ``decoupled
ODEs" by changing $\alpha$ and $\beta$ into new variables $y$ and
$q$, where \beq\left.\begin{array}{l}y=
m^{-\frac{3(3-\gamma)}{2(3\gamma-1)}}z^{\frac{\gamma+1}{2(\gamma-1)}}
((u+mz)_x -\frac{2}{3\gamma-1}m_x z),\label{intr
main1}\end{array}\right.\eeq \beq\left.\begin{array}{l}q=
m^{-\frac{3(3-\gamma)}{2(3\gamma-1)}}z^{\frac{\gamma+1}{2(\gamma-1)}}
((u-mz)_x +\frac{2}{3\gamma-1}m_x z). \label{intr
main2}\end{array}\right.\eeq Similarly, we define
\beq\left.\begin{array}{l}\tilde{y}=
z^{\frac{\gamma+1}{2(\gamma-1)}}((u+m z)_x -\frac{2}{3\gamma-1}m_x
z),\label{intr main3}\end{array}\right.\eeq
\beq\left.\begin{array}{l}\tilde{q}=
z^{\frac{\gamma+1}{2(\gamma-1)}}((u-m z)_x +\frac{2}{3\gamma-1}m_x
z).\label{intr main4}\end{array}\right.\eeq {\theorem{\label{them
new ODEs}The smooth solutions in
(\ref{lagrangian1})}$\sim$(\ref{lagrangian3}) satisfy
\begin{eqnarray} y^\prime&=&a_0+ a_2 y^2, \label{new ode1}
\\q^\backprime&=&a_0+ a_2 q^2, \label{new ode2}\end{eqnarray}
where
\begin{eqnarray}{a}_0=&\frac{K_c}{\gamma}
m^{-\frac{3(3-\gamma)}{2(3\gamma-1)}} [\frac{\gamma-1}{3\gamma-1}m
m_{xx}-\frac{(3\gamma+1) (\gamma-1)}{(3\gamma-1)^2}m_x^2]
z^{\frac{3(\gamma+1)}{2(\gamma-1)}+1},\label{a0def}\\
{a}_2=&-K_c\frac{\gamma+1}{2(\gamma-1)}
m^{\frac{3(3-\gamma)}{2(3\gamma-1)}}
z^{\frac{\gamma+1}{2(\gamma-1)}-1}<0.\qquad\qquad\qquad\qquad\label{a2def}\end{eqnarray}
}
\\
Clearly, the functions $a_0$ and $a_2$ only depend on the density
and initial entropy profile. These ``ODEs" generalize {Lax's}
``ODEs" for $p$-system in \cite{lax2}. In fact, when entropy is
constant, $a_0=0$, and we recover Lax's equations. In this paper,
``$ \rightarrow $" denotes the limit when $(t,x)$ approaches a fixed
(finite) point $(t_0,x_0)$.

For later reference, we state two assumptions that will be used in
some of our results.
\begin{eqnarray} Assumption\ 1:z\ is\ not\ equal\ to\ zero\ or\ infinity\
for\ any \ x\ and\
t.\qquad\qquad\qquad\nonumber\end{eqnarray}\begin{eqnarray}
Assumption\ 2: There\ exist\ positive\ constants\ Z_L,\ Z_U,\ M_1
\sim M_4,\ such\ that\qquad\qquad\nonumber\end{eqnarray}
\begin{eqnarray}
Z_L<z<Z_U,\ for\ all\ x\ and\ t;\label{assumption2Z}
\end{eqnarray}
\begin{eqnarray}
M_1<m^0<M_2,\quad |m^0_x|<M_3,\quad
|m^0_{xx}|<M_4,\label{assumption2m}
\end{eqnarray}  \emph{where
$m^0(x)=m(0,x)$ is given by the prescribed entropy profile
$s(0,x)\equiv s(t,x)$ by (\ref{s con}) and (\ref{intr z
m}).}\vspace{.2cm}
\\
By (\ref{intr z m}), Assumption 1 means that the density is not zero
or infinity, and (\ref{assumption2Z}) means that the density has
positive upper and lower bounds. Furthermore, all realistic smooth
($C^2$) initial entropy profiles satisfy (\ref{assumption2m}). We
impose $Z_L$ to avoid potential problems at vacuum, which are
addressed in an upcoming paper \cite{young com}. $Z_L$ can be
arbitrarily small.

Using (\ref{new ode1}) and (\ref{new ode2}), we can now give some
singularity formation results. {\theorem{\label{Thm singularity2}
Assume the initial data are smooth, and Assumption 2 holds. The
constants $Z_U$ and $M_1$$\sim$$M_4$ are given in Assumption 2. Then
there exist positive constants $N$ and $\tilde{N}$ depending only on
$Z_U$ and $M_1$$\sim$$M_4$, such that, if $y$ or $q$ is less than
$-N$ or $\tilde{y}$ or $\tilde{q}$ is less than $-\tilde{N}$
somewhere in the initial data, then $|u_x|$ and/or $|\tau_x|$ blow
up in finite time. When $\gamma\geqslant3$, all the results hold
without the lower bound of density in Assumption 2. (Here $y$, $q$,
$\tilde{y}$ and $\tilde{q}$ are defined in (\ref{intr
main1})$\sim$(\ref{intr main4}).)}} {\theorem{\label{singularity}
Assume the initial data are smooth, and Assumption 2 holds.
Furthermore, assume that there exists a point $A\in \mathbb{R}$,
such that the initial entropy profile satisfies
\beq(m^{-\frac{2}{3\gamma-1}})_{xx}\geqslant 0,\ for \
x>A.\label{condition1}\eeq  If, in the initial data, $y_0=y(0,
A^*)<0$ at some point $x=A^*$ with $A^*>A$, then $|u_x|$ and/or
$|\tau_x|$ blow up before some finite time $T_*$, where \beq
0<T_*\leqslant-\frac{1}{y_0 \min(-a_2)},\eeq and $\min(-a_2)$ is a
positive constant depending on $Z_L\ (Z_U)$ and $M_1\ (M_2)$ when
$1<\gamma<3\ (\gamma>3)$. When $\gamma=3$, all the results hold
without Assumption 2. Symmetric results hold for $q$. (Here $a_2$ is
defined in (\ref{a2def}).)}}

Theorem \ref{Thm singularity2} implies that gradients of solutions
blow up if the initial compressions are strong enough, which is a
direct generalization of Lax's singularity formation results for
$p$-system. When the variation of entropy is mild, $N$ and
$\tilde{N}$ are close to zero, so the shock free solutions are
``almost rarefactive", which is consistent with Lax's singularity
formation results.   In \cite{young blake 1} and a forthcoming paper
\cite{G1}, examples are given, of solutions containing compressive
waves, but the gradients of the solutions don't blow up.

Our ``ODEs" and the singularity formation results are for
arbitrarily large smooth initial data, where there is no restriction
on the amplitude of the waves.

We divide this paper into $7$ sections. In sections $2$ and $3$, we
review the background and definition of $R/C$ character in the
$p$-system. In section $4$, we define rarefactive and compressive
waves in a smoothly varying entropy field. In sections $5$ and $6$,
we give ``ODEs" for smooth solutions. In section $7$, we prove the
singularity formation results.
%
%
\section{Equations and Coordinates}
In this paper, we focus on the polytropic ideal gas dynamics, where
the equation of state is given by, \beq p\tau=RT,\label{introduction
1}\eeq and \beq e=c_v
T={\frac{1}{\gamma-1}}p\tau,\label{introduction 2}\eeq  with \beq
p=Ke^{\frac{S}{c_v}}\tau^{-\gamma}.\label{introduction 3}\eeq Here
$S$ is the entropy, $T$ is the temperature, $R$, $K$, $c_v$ are all
positive constants, and the adiabatic gas constant $\gamma>1$, c.f.
\cite{courant}. The (Lagrangian) wave (sound) speed is given by \beq
c=\sqrt{-p_\tau}=\sqrt{K\gamma}{\tau}^{-\frac{\gamma+1}{2}}e^
{\frac{S}{2c_v}}.\label{c def}\eeq

We use the coordinates provided by {B. Temple} and {R. Young}, in
\cite{young blake 1}. Define new variables $m$ and $z$ for $S$ and
$\tau$, by \beq m=e^{\frac{S}{2c_v}},\label{m def}\eeq and \beq
z=\int^\infty_\tau{\frac{c}{m}d\tau}=\frac{2\sqrt{K\gamma}}{\gamma-1}
\tau^{-\frac{\gamma-1}{2}},\label{z def}\eeq where we use (\ref{c
def}) and (\ref{m def}).  It follows that
\begin{eqnarray}\tau&=&K_{\tau}z^{-\frac{2}{\gamma-1}},\label{tau p c 1}\\
p&=&K_p m^2 z^{\frac{2\gamma}{\gamma-1}},\label{tau p c 2}\\
c&=&c(z,m)=K_c m z^{\frac{\gamma+1}{\gamma-1}},
 \label{tau p c 3}\end{eqnarray} where $K_\tau$, $K_p$ and $K_c$ are positive
constants given by
\begin{eqnarray}K_\tau&=&(\frac{2\sqrt{K\gamma}}{\gamma-1})
^\frac{2}{\gamma-1},\label{constant1}\\ K_p&=&K
K_\tau^{-\gamma},\label{constant2}\\
K_c&=&\sqrt{K\gamma}K_\tau^{-\frac{\gamma+1}{2}},\label{constant3}\end{eqnarray}
and \beq K_p=\frac{\gamma-1}{2\gamma}K_c.\label{KpKcRela}\eeq By
(\ref{m def}) and (\ref{z def}), for $C^1$ solutions, the
{Lagrangian} equations (\ref{lagrangian1})$\sim$(\ref{lagrangian3})
are equivalent to
\begin{eqnarray}z_t+\frac{c}{m}u_x=0,\label{lagrangian1 zm}
\\u_t+mcz_x+2\frac{p}{m}m_x=0,\label{lagrangian2 zm}\\
m_t=0\label{lagrangian3 zm},
\end{eqnarray}
where the last equation is coming from (\ref{s con}) instead of
(\ref{lagrangian3}), c.f. \cite{smoller}.

When the entropy is constant, $p=p(\tau)$, so (\ref{lagrangian1
zm})$\sim$(\ref{lagrangian3 zm}) change to the $p$-system:
\begin{eqnarray}z_t+\frac{c}{m}u_x=0,\label{p1 zm}
\\u_t+mcz_x=0.\label{p2 zm}
\end{eqnarray}
So the corresponding Riemann invariants are
\begin{eqnarray}r=u-mz,\label{r def}\\s=u+mz,\label{s def}\end{eqnarray}
which satisfy, by (\ref{p1 zm}) and (\ref{p2 zm}),
\begin{eqnarray} r_t-cr_x=0,\label{r invariant}\\
s_t+cs_x=0.\label{s invariant}\end{eqnarray}
%
%
\section{Rarefactive and Compressive Waves in a Constant Entropy Domain}
We first review the definition of the rarefactive and compressive
waves for smooth solutions in a constant entropy domain, where the
system is $2\times 2$. At this time, the system has full set of
Riemann coordinates (\ref{r def}) and (\ref{s def}). For a single
simple wave, the wave is rarefactive (compressive) if pressure,
density or wave speed, are decreasing (increasing) from ahead of to
behind the wave, c.f. \cite{lax}. Here we consider not only the
simple waves, but also the local rarefactive and compressive
characters at some points or open sets in the wave interaction
regions, c.f. \cite{lax2}, \cite{young blake 1}, \cite{young Global
interaction}. Sometimes, we use the names ``rarefactive wave" and
``compressive waves" or even ``\emph{R}" and ``\emph{C}", instead of
the rarefactive and compressive character at some points or open
sets, without confusion.

In this section, we assume entropy is a constant, hence $p$ is a
function which only depends on $\tau$, so that \beq
p^{\backprime}=p_\tau \tau^{\backprime},\quad c^{\backprime}=c_\tau
\tau^{\backprime}, \quad
\rho^{\backprime}=-\frac{1}{\tau^2}\tau^{\backprime}.\label{p system
judge1}\eeq Recall subscripts $``\prime"$ and $``\backprime"$ denote
the directional derivatives along forward and backward
characteristics, respectively, which are defined in (\ref{prime
backprime}). Note $p_\tau<0$, and $c_\tau<0$ since $p_{\tau\tau}>0$
and (\ref{c def}). Furthermore, \beq -\frac{s_t}{c}=s_x=u_x-c
\tau_x=\tau^{\backprime},\label{p system judge2}\eeq where we use
(\ref{s invariant}), (\ref{s def}) and (\ref{lagrangian1}),
respectively.

Hence \beq \tau^{\backprime}<0\Leftrightarrow
c^{\backprime}>0\Leftrightarrow p^{\backprime}>0 \Leftrightarrow
s_x<0\Leftrightarrow s_t>0\Leftrightarrow
\rho^{\backprime}>0,\label{simpleRC 1}\eeq which means that the
forward waves are compressive; and \beq
\tau^{\backprime}>0\Leftrightarrow c^{\backprime}<0\Leftrightarrow
p^{\backprime}<0 \Leftrightarrow s_x>0\Leftrightarrow
s_t<0\Leftrightarrow \rho^{\backprime}<0,\label{simpleRC 2}\eeq
which means that the forward waves are rarefactive.

In fact, (\ref{simpleRC 1}) and (\ref{simpleRC 2}) mean that
pressure, density and wave speed, are locally increasing or
decreasing from ahead of to behind the forward wave along backward
characteristic line, respectively.  We use the directional
derivative along backward characteristic to define the forward
\emph{R/C} character, which discounts the disturbance from the
backward waves. Symmetrically, we define the backward \emph{R/C}
character.
%
%
\section{Rarefactive and Compressive Solutions in a Smoothly Varying
Entropy Domain} In the previous Section, we describe the equivalent
definitions of rarefactive and compressive waves in a constant
entropy domain. However, when the entropy is smoothly varying, there
are no Riemann invariants anymore, and the forward and backward
waves can't be geometrically divided, which give difficulties for
defining the R/C characters. In this
 Section, we present a definition of rarefactive and compressive
waves in a smoothly varying entropy domain, which is new and
fundamental for the smooth waves of the compressible Euler
equations.

The conditions shown in the previous Section are not equivalent
anymore in the domain with smoothly varying entropy profile. This is
because pressure and wave speed become functions depending on both
density and entropy.  In order to pick up the right condition, we
first consider the stationary solutions of
(\ref{lagrangian1})$\sim$(\ref{lagrangian3}), where $u$, $\tau$, $S$
are stationary. So $u$, $p$ are constant by
(\ref{lagrangian1})$\sim$(\ref{lagrangian3}). In particular, the
density need not be constant when the entropy profile is varying.

We go back to check the rarefactive and compressive conditions given
in the previous Section. In the stationary solutions,
$p^\prime=p^\backprime=0$ while the other quantities are nonzero if
the entropy profile is smoothly varying. In the other words,
pressure is the only thermodynamic variable which is not impacted by
the stationary entropy field. Moreover, change of pressure
contributes to the wave propagation on the (non-vertical)
characteristic lines, i.e. rarefactive or compressive waves. So
pressure is the only possible variable which still can distinguish
rarefactive and compressive waves in a smoothly varying entropy
field. So we can get the definition of rarefactive and compressive
waves in Definition \ref{def1}. Note that we use the directional
derivative of pressure along the \emph{opposite} characteristic to
define the \emph{R/C} character, which helps us discount the
disturbance from waves in the opposite direction. This definition
also applies for the general pressure law with $p_\tau<0,\
p_{\tau\tau}>0$.

We define \beq \alpha=-\frac{p^\backprime}{c^2},\quad
\beta=-\frac{p^\prime}{c^2}.\label{alpha beta p}\eeq {\lemma
\label{lemma alpha beta p}\beq \left.\begin{array}{l}
\alpha=u_x+mz_x+\frac{\gamma-1}{\gamma}m_x
z,\end{array}\right.\label{alpha def}\eeq and \beq
\left.\begin{array}{l}\beta=u_x-mz_x-\frac{\gamma-1}{\gamma}m_x
z.\end{array}\right.\label{beta def}\eeq }

\begin{proof}
\beq\left.\begin{array}{lll} -c^2 \alpha&=&p^\backprime\\&=&(K_p m^2
z^{\frac{2\gamma}{\gamma-1}})^\backprime\\
&=&K_p m^2 \frac{2\gamma}{\gamma-1}
z^{\frac{\gamma+1}{\gamma-1}}z_t-c K_p m^2 \frac{2\gamma}{\gamma-1}
z^{\frac{\gamma+1}{\gamma-1}}z_x-2c K_p m m_x
z^{\frac{2\gamma}{\gamma-1}}\\ &=& -K_p m c \frac{2\gamma}{\gamma-1}
z^{\frac{\gamma+1}{\gamma-1}}(u_x+mz_x+\frac{\gamma-1}{\gamma}m_x
z)\\
&=&c^2(u_x+mz_x+\frac{\gamma-1}{\gamma}m_x z),
\end{array}\right.\label{p alpha beta}\eeq
where we use (\ref{tau p c 3}) and (\ref{KpKcRela}). Similarly, we
can prove (\ref{beta def}).
\end{proof}

By the Definition \ref{def1} and (\ref{alpha beta p}), we can
equivalently define the $R/C$ character by $\alpha$ and $\beta$.
{\lemma{\label{def2}}{In the polytropic ideal gas, the local {R/C}
character of the smooth solution is given by: \beq
\left.\begin{array}{llll} \it{Forward}& R &
\it{iff}& \alpha>0,\\ \it{Forward}& C & \it{iff} & \alpha<0,\\
\it{Backward}& R & \it{iff}& \beta>0,\\ \it{Backward}& C &\it{iff}&
\beta<0.
\end{array}\right.\eeq
When Assumption 1 holds, \beq\left.\begin{array}{l} |\alpha|\
\it{or}\ |\beta|\rightarrow \infty \ \it{iff}\ |u_x|\ \it{or}\
|\tau_x| \rightarrow \infty.\label{alpha beta
blowup}\end{array}\right.\eeq}}
\begin{proof}
Clearly \beq\left.\begin{array}{l} p^\backprime \gtrless 0
\Leftrightarrow \alpha \lessgtr 0,\end{array}\right.\eeq and,
\beq\left.\begin{array}{l} p^\prime \gtrless 0 \Leftrightarrow \beta
\lessgtr 0.\end{array}\right.\eeq

By (\ref{alpha def}) and
(\ref{beta def}), \begin{eqnarray} \alpha+\beta&=&2u_x,\label{alphabetaux}\\
\alpha-\beta&=&2(mz_x+\frac{\gamma-1}{\gamma}m_x
z).\label{alphabezxmx}\end{eqnarray} By (\ref{alpha def}),
(\ref{beta def}), (\ref{alphabetaux}), (\ref{alphabezxmx}), (\ref{z
def}) and Assumption 1, we get (\ref{alpha beta blowup}), where we
also use that $m$ is stationary and positive in smooth solution by
(\ref{m def}) and (\ref{lagrangian3 zm}).
\end{proof}

Note: We give the definition of compressive and rarefactive waves in
a smoothly varying entropy profile from a purely physical point of
view, which also explains the physical meanings of Riemann
invariants in the $p$-system. In fact, in the $p$-system, the
physical meanings of the derivatives of Riemann invariants can be
explained by the directional derivatives of pressure: \beq
s_x=-\frac{p^\backprime}{c^2},\quad
r_x=-\frac{p^\prime}{c^2},\label{Riemann physical}\eeq where we use
that $m$ is a constant, (\ref{r def}), (\ref{s def}), (\ref{alpha
beta p}) and Lemma \ref{lemma alpha beta p}. So $\alpha$, $\beta$
can be considered as the generalization of $s_x$ and $r_x$ in a
smoothly varying entropy field. When $m_x=0$, $\alpha$ and $\beta$
equal to $s_x$ and $r_x$, respectively.
%
%
\section{Smooth Wave Propagation}
In this section, we consider the smooth wave propagation of
(\ref{lagrangian1})$\sim$(\ref{lagrangian3}). By considering the
directional derivatives of $\alpha$ and $\beta$, we construct
``ordinary differential equations" (\ref{rem1}) and (\ref{rem2}) in
Theorem \ref{remark}, which together with the definition of
rarefactive and compressive waves will give us a framework for the
wave propagation of smooth solutions. These ``ODEs" are not real
ODEs since the derivatives are in different directions and
(\ref{rem1}) and (\ref{rem2}) depend on the variable $z$, but they
are much simpler than the original partial differential equations.
\subsection{Proof of Theorem \ref{remark}}
\begin{proof}
By Lemma \ref{lemma alpha beta p}, \beq \left.\begin{array}{ll}
\alpha^\prime&=(u_x+mz_x+\frac{\gamma-1}{\gamma}m_x
z)^{\prime}\\&=u_{xt}+m z_{xt}+\frac{\gamma-1}{\gamma}m_x z_t+c[
u_{xx}+ m z_{xx}+ m_{x}z_{x}+\frac{\gamma-1}{\gamma}( m_x z_{x}+
m_{xx}z)]\\&=(u_{xt}+c m z_{xx})+(c u_{xx}+m
z_{xt})+\frac{2\gamma-1}{\gamma}c m_x z_x+\frac{\gamma-1}{\gamma}m_x
z_t+\frac{\gamma-1}{\gamma}c
m_{xx}z.\end{array}\right.\label{app0}\eeq By (\ref{lagrangian1
zm})$\sim$(\ref{lagrangian3 zm}) and (\ref{tau p c 2}),
\beq(u_t+cmz_x)_x=(-2\frac{p}{m}m_x)_x=-(2K_p m m_x
z^{\frac{2\gamma}{\gamma-1}})_x,\eeq and \beq (cu_x+m z_t)_x=0.\eeq
So \beq \left.\begin{array}{ll}u_{xt}+cmz_{xx}&=-2K_c m m_x
z^{\frac{\gamma+1}{\gamma-1}}z_x-\frac{\gamma+1}{\gamma-1}K_c m^2
z^{\frac{2}{\gamma-1}}(z_x)^2\\&-2K_p ({m_x})^2
z^{\frac{2\gamma}{\gamma-1}}-2K_p m
m_{xx}z^{\frac{2\gamma}{\gamma-1}}-{\frac{4\gamma}{\gamma-1}}K_p m
m_{x}z^{\frac{\gamma+1}{\gamma-1}}z_x,\end{array}\right.\label{app1}\eeq
and \beq \left.\begin{array}{l}  cu_{xx}+m z_{xt}=-c_x u_x-m_x
z_t=-{\frac{\gamma+1}{\gamma-1}}K_c m z^{\frac{2}{\gamma-1}}z_x
u_x.\end{array}\right.\label{app2}\eeq Furthermore, the sum of the
second term in the right hand side of (\ref{app1}) and the right
hand side of (\ref{app2}) is \beq \left.\begin{array}{ll}
&-\frac{\gamma+1}{\gamma-1}K_c m z^{\frac{2}{\gamma-1}}z_x(m
z_x+u_x)\\=&-\frac{\gamma+1}{\gamma-1}K_c
z^{\frac{2}{\gamma-1}}(\alpha-\frac{\gamma-1}{\gamma}m_x z-u_x)
(\alpha-\frac{\gamma-1}{\gamma}m_x z)\\=&-\frac{\gamma+1}{\gamma}K_c
m_x z^{\frac{\gamma+1}{\gamma-1}}u_x-\frac{\gamma^2-1}{\gamma^2}K_c
(m_x)^2 z^\frac{2\gamma}{\gamma-1}
\\& +\frac{\gamma+1}{\gamma-1}K_c z^{\frac{2}{\gamma-1}}
[(\frac{2\gamma-2}{\gamma}m_x z+u_x)\alpha-\alpha^2].
\end{array}\right.\label{app3}\eeq

By (\ref{app1})$\sim$(\ref{app3}) and (\ref{tau p c 3}), the right
hand side of (\ref{app0}) equals to, \beq \left.\begin{array}{l}
-2K_c m m_x z^{\frac{\gamma+1}{\gamma-1}}z_x-2K_p (m_{x})^2
z^{\frac{2\gamma}{\gamma-1}}-2K_p m
m_{xx}z^{\frac{2\gamma}{\gamma-1}}-{\frac{4\gamma}{\gamma-1}}K_p m
m_{x}z^{\frac{\gamma+1}{\gamma-1}}z_x\\-\frac{\gamma+1}{\gamma}K_c
m_x z^{\frac{\gamma+1}{\gamma-1}}u_x-\frac{\gamma^2-1}{\gamma^2}K_c
(m_x)^2 z^\frac{2\gamma}{\gamma-1}\\+\frac{2\gamma-1}{\gamma}K_c m
m_x z^\frac{\gamma+1}{\gamma-1} z_x-\frac{\gamma-1}{\gamma}K_c m_x
z^\frac{\gamma+1}{\gamma-1}u_x+\frac{\gamma-1}{\gamma}K_c m m_{xx}
z^{\frac{2\gamma}{\gamma-1}} \\
+\frac{\gamma+1}{\gamma-1}K_c z^{\frac{2}{\gamma-1}}
[(\frac{2\gamma-2}{\gamma}m_x z+u_x)\alpha-\alpha^2],
\end{array}\right.\label{app4}\eeq
where we use (\ref{lagrangian1 zm}) to get rid of $z_t$. Plug
(\ref{KpKcRela}) into (\ref{app4}), then (\ref{app4}) equals to,\beq
\left.\begin{array}{l} -2K_c m_x
z^{\frac{\gamma+1}{\gamma-1}}u_x-\frac{2\gamma^2-\gamma-1}{\gamma^2}K_c
(m_x)^2 z^\frac{2\gamma}{\gamma-1}+\frac{-2\gamma-1}{\gamma}K_c m
m_x z^\frac{\gamma+1}{\gamma-1} z_x \\
+\frac{\gamma+1}{\gamma-1}K_c z^{\frac{2}{\gamma-1}}
[(\frac{2\gamma-2}{\gamma}m_x z+u_x)\alpha-\alpha^2].
\end{array}\right.\label{app5}\eeq
Using $ m z_x =\alpha-u_x-\frac{\gamma-1}{\gamma}m_x z,$ which is
from (\ref{alpha def}), (\ref{app5}) can be simplified to
\beq\left.\begin{array}{l}K_c
z^{\frac{\gamma+1}{\gamma-1}}\{\frac{1}{\gamma}m_x u_x+
(\frac{1}{\gamma}m_x+\frac{\gamma+1}{\gamma-1}\frac{u_x}{z})\alpha-
\frac{\gamma+1}{\gamma-1}\frac{\alpha^2}{z}\}.\end{array}\right.\label{app6}\eeq
By (\ref{alphabetaux}), \beq
u_x=\frac{\alpha+\beta}{2}.\label{uxalphabetahalf}\eeq Plugging
(\ref{uxalphabetahalf}) into (\ref{app6}),  we get (\ref{rem1}). The
calculation of (\ref{rem2}) is same.
\end{proof}

{\cor\label{them5} {In an open set $U$ in the $(t,x)$-plane, assume
the solutions of (\ref{lagrangian1})$\sim$(\ref{lagrangian3}) are
smooth, $m_x\neq 0$, and Assumption 1 holds. If $\beta\equiv 0$ in
$U$, then $\alpha\equiv 0$ in $U$; if $\alpha\equiv 0$ in $U$, then
$\beta\equiv 0$ in $U$.}}\vspace{.1cm}\begin{proof} If $\beta\equiv
0$ in $U$, $-k_1 k_2\alpha\equiv 0$ by (\ref{rem2}), hence
$\alpha\equiv 0$ by Assumption 1. Another case can be proved
similarly.
\end{proof}

In this corollary, we show that there will be no ``pure" forward or
backward rarefactive or compressive waves inside any open set in
$(t,x)$-plane with varying entropy, i.e. the forward and backward
waves will exist or die out mutually. This is different from the
isentropic domain, where we can find rarefactive or compressive
simple waves.
%
%
\section{``Decoupled" ODEs}
Next, after introducing new variables $y$ and $q$, we change
(\ref{rem1}) and (\ref{rem2}) into ``decoupled ODEs" (\ref{new
ode1}) and (\ref{new ode2}) in Theorem \ref{Thm singularity2}, which
are only coupled by $z$ and $x$. The variables $y$ and $q$ are given
in (\ref{intr main1}) and (\ref{intr main2}). These new ``ODEs" will
help us prove the singularity results. The equations (\ref{new
ode1}) and (\ref{new ode2}) generalize {Lax's} ``ODEs" for
$p$-system in \cite{lax2}. In fact, when the entropy is constant,
$a_0=0$, so (\ref{new ode1}) and (\ref{new ode2}) change to \beq
y^\prime= a_2 y^2,\quad q^\backprime= a_2 q^2,\eeq which is exactly
the ``ODEs" provided in \cite{lax2} for the $p$-system.
\subsection{Proof of Theorem \ref{them
new ODEs}} \begin{proof} By (\ref{lagrangian1 zm}), (\ref{tau p c
3}) and (\ref{beta def}),
\begin{eqnarray}z^{\prime}&=&z_t+c z_x\nonumber\\
&=&-\frac{c}{m}(u_x-m z_x)\nonumber\\&=&-K_c
z^{\frac{\gamma+1}{\gamma-1}}(\beta+\frac{\gamma-1}{\gamma}m_x z).
\end{eqnarray} Hence \beq \beta=-\frac{z^{\prime}}
{K_c z^{\frac{\gamma+1}{\gamma-1}}}-\frac{\gamma-1}{\gamma}m_x
z.\label{beta z prime}\eeq Plugging (\ref{beta z prime}) into
(\ref{rem1}), we get
\beq\left.\begin{array}{l}\alpha^\prime=k_1\{k_2
(3\alpha-\frac{z^{\prime}} {K_c
z^{\frac{\gamma+1}{\gamma-1}}}-\frac{\gamma-1}{\gamma}m_x
z)+\alpha(-\frac{z^{\prime}} {K_c
z^{\frac{\gamma+1}{\gamma-1}}}-\frac{\gamma-1}{\gamma}m_x
z)-\alpha^2\}.\end{array}\right.\eeq We move the terms including
$z^\prime$ to the left hand side, then we multiply by
$z^{\frac{\gamma+1}{2(\gamma-1)}}$ on both sides. After
simplification, we have \beq\left.\begin{array}{ll} &
z^{\frac{\gamma+1}{2(\gamma-1)}}\alpha^{\prime}+
\frac{m_x}{2\gamma}z^{\frac{\gamma+1}{2(\gamma-1)}}z^{\prime}+
\frac{\gamma+1}{2(\gamma-1)}{\alpha}
z^{\frac{\gamma+1}{2(\gamma-1)}-1}z^{\prime}\\=& -\frac{\gamma-1}{2
\gamma^2}K_c z^{\frac{3(\gamma+1)}{2(\gamma-1)}+1}m_x^2
+\frac{2-\gamma}{2 \gamma}K_c m_x
z^{\frac{3(\gamma+1)}{2(\gamma-1)}}\alpha-
\frac{\gamma+1}{2(\gamma-1)}K_c
z^{\frac{3(\gamma+1)}{2(\gamma-1)}-1}
\alpha^2,\end{array}\right.\label{y eqn 1}\eeq where we use (\ref{k
def}). The left hand side of (\ref{y eqn 1}) is equal to
\beq\left.\begin{array}{l} (\alpha z^{\frac{\gamma+1}{2(\gamma-1)}}+
\frac{\gamma-1}{\gamma(3\gamma-1)}
z^{\frac{\gamma+1}{2(\gamma-1)}+1}
m_x)^{\prime}-\frac{\gamma-1}{\gamma(3\gamma-1)}z^{\frac{\gamma+1}{2(\gamma-1)}+1}
(m_x)^\prime.\label{y eqn 2}\end{array}\right.\eeq We define a new
variable $\tilde{y}, $\beq \tilde{y}=\alpha
z^{\frac{\gamma+1}{2(\gamma-1)}}+ \frac{\gamma-1}{\gamma(3\gamma-1)}
z^{\frac{\gamma+1}{2(\gamma-1)}+1} m_x.\label{def tilde y}\eeq So
\beq \alpha=\tilde{y}
z^{-\frac{\gamma+1}{2(\gamma-1)}}-\frac{\gamma-1}{\gamma(3\gamma-1)}
z m_x. \label{alpha tilde y 1}\eeq Hence, by (\ref{y eqn 2}),
(\ref{alpha tilde y 1}) and $(m_x)^\prime=c m_{xx}$, (\ref{y eqn 1})
changes to
\begin{eqnarray}\tilde{y}^\prime=&\frac{\gamma-1}{\gamma(3\gamma-1)}K_c
z^{\frac{3(\gamma+1)}{2(\gamma-1)}+1}m
m_{xx}-\frac{\gamma-1}{2\gamma^2}K_c
z^{\frac{3(\gamma+1)}{2(\gamma-1)}+1}m_x^2\nonumber\\&+\frac{2-\gamma}{2
\gamma}K_c m_x z^{\frac{3(\gamma+1)}{2(\gamma-1)}}(\tilde{y}
z^{-\frac{\gamma+1}{2(\gamma-1)}}-\frac{\gamma-1}{\gamma(3\gamma-1)}
z m_x)\nonumber\\&- \frac{\gamma+1}{2(\gamma-1)}K_c
z^{\frac{3(\gamma+1)}{2(\gamma-1)}-1} (\tilde{y}
z^{-\frac{\gamma+1}{2(\gamma-1)}}-\frac{\gamma-1}{\gamma(3\gamma-1)}
z m_x)^2.\label{y eqn 3}
\end{eqnarray}
 After simplification,
we have \beq \tilde{y}^\prime=\tilde{a}_0 +\tilde{a}_1 \tilde{y}+
\tilde{a}_2 \tilde{y}^2,\label{y fun}\eeq where
\begin{eqnarray}\tilde{a}_0&=&K_c\frac{1}{\gamma}[\frac{\gamma-1}
{3\gamma-1}m m_{xx}-\frac{(3\gamma+1)
(\gamma-1)}{(3\gamma-1)^2}m_x^2]
z^{\frac{3(\gamma+1)}{2(\gamma-1)}+1},\label{tilde a0}\\
\tilde{a}_1&=&K_c\frac{3(3-\gamma)}{2(3\gamma-1)}
m_x z^{\frac{\gamma+1}{\gamma-1}},\label{tilde a1}\\
\tilde{a}_2&=&-K_c\frac{\gamma+1}{2(\gamma-1)}
z^{\frac{\gamma+1}{2(\gamma-1)}-1}.\label{tilde a2}\end{eqnarray}
Then we do one more simplification by multiplying
\beq{\bar{\mu}}=m^{-\frac{3(3-\gamma)}{2(3\gamma-1)}},\label{mu}\eeq
on (\ref{y fun}). In fact, it is easy to check that \beq
{\bar{\mu}}^{\prime}=-\tilde{a}_1 {\bar{\mu}},\eeq since $
m^{\prime}=cm_x.$ Then we denote \beq
y={\bar{\mu}}\tilde{y}\label{def y}.\eeq Hence (\ref{y fun}) changes
to \beq y^\prime=a_0+a_2 y^2,\eeq where \beq
a_0={\bar{\mu}}\tilde{a}_0,\
a_2=\tilde{a}_2/{\bar{\mu}}.\label{a0a2}\eeq\vspace{.1cm}

Similarly, we prove (\ref{new ode2}). By (\ref{lagrangian1 zm}),
(\ref{tau p c 3}) and (\ref{beta def}),
\begin{eqnarray}z^{\backprime}&=&z_t-c z_x\nonumber\\
&=&-\frac{c}{m}(u_x+m z_x)\nonumber\\&=&-K_c
z^{\frac{\gamma+1}{\gamma-1}}(\alpha-\frac{\gamma-1}{\gamma}m_x z).
\end{eqnarray} Hence \beq \alpha=-\frac{z^{\backprime}}
{K_c z^{\frac{\gamma+1}{\gamma-1}}}+\frac{\gamma-1}{\gamma}m_x
z.\label{alpha z backprime}\eeq Plugging (\ref{alpha z backprime})
into (\ref{rem2}), we get
\beq\left.\begin{array}{l}\beta^\backprime=k_1\{-k_2
(3\beta-\frac{z^{\backprime}} {K_c
z^{\frac{\gamma+1}{\gamma-1}}}+\frac{\gamma-1}{\gamma}m_x
z)+\beta(-\frac{z^{\backprime}} {K_c
z^{\frac{\gamma+1}{\gamma-1}}}+\frac{\gamma-1}{\gamma}m_x
z)-\beta^2\}.\end{array}\right.\eeq We move the terms including
$z^\backprime$ to the left hand side, then we multiply by
$z^{\frac{\gamma+1}{2(\gamma-1)}}$ on both sides. After
simplification, we have \beq\left.\begin{array}{ll} &
z^{\frac{\gamma+1}{2(\gamma-1)}}\beta^{\backprime}-
\frac{m_x}{2\gamma}z^{\frac{\gamma+1}{2(\gamma-1)}}z^{\backprime}+
\frac{\gamma+1}{2(\gamma-1)}{\beta}
z^{\frac{\gamma+1}{2(\gamma-1)}-1}z^{\backprime}\\=&
-\frac{\gamma-1}{2 \gamma^2}K_c
z^{\frac{3(\gamma+1)}{2(\gamma-1)}+1}m_x^2 -\frac{2-\gamma}{2
\gamma}K_c m_x z^{\frac{3(\gamma+1)}{2(\gamma-1)}}\beta-
\frac{\gamma+1}{2(\gamma-1)}K_c
z^{\frac{3(\gamma+1)}{2(\gamma-1)}-1}
\beta^2,\end{array}\right.\label{q eqn 1}\eeq where we use (\ref{k
def}). The left hand side of (\ref{q eqn 1}) is equal to
\beq\left.\begin{array}{l} (\beta z^{\frac{\gamma+1}{2(\gamma-1)}}-
\frac{\gamma-1}{\gamma(3\gamma-1)}
z^{\frac{\gamma+1}{2(\gamma-1)}+1}
m_x)^{\backprime}+\frac{\gamma-1}{\gamma(3\gamma-1)}z^{\frac{\gamma+1}
{2(\gamma-1)}+1}(m_x)^\backprime.\label{q eqn
2}\end{array}\right.\eeq We define a new variable $\tilde{q}$, \beq
\tilde{q}=\beta z^{\frac{\gamma+1}{2(\gamma-1)}}-
\frac{\gamma-1}{\gamma(3\gamma-1)}
z^{\frac{\gamma+1}{2(\gamma-1)}+1} m_x.\label{def tilde q}\eeq So
\beq \beta=\tilde{q}
z^{-\frac{\gamma+1}{2(\gamma-1)}}+\frac{\gamma-1}{\gamma(3\gamma-1)}
z m_x.\eeq Hence (\ref{q eqn 1}) changes to
\begin{eqnarray}\tilde{q}^\backprime=&\frac{\gamma-1}{\gamma(3\gamma-1)}K_c
z^{\frac{3(\gamma+1)}{2(\gamma-1)}+1}m
m_{xx}-\frac{\gamma-1}{2\gamma^2}K_c
z^{\frac{3(\gamma+1)}{2(\gamma-1)}+1}m_x^2\nonumber\\&-\frac{2-\gamma}{2
\gamma}K_c m_x z^{\frac{3(\gamma+1)}{2(\gamma-1)}}(\tilde{q}
z^{-\frac{\gamma+1}{2(\gamma-1)}}+\frac{\gamma-1}{\gamma(3\gamma-1)}
z m_x)\nonumber\\&- \frac{\gamma+1}{2(\gamma-1)}K_c
z^{\frac{3(\gamma+1)}{2(\gamma-1)}-1} (\tilde{q}
z^{-\frac{\gamma+1}{2(\gamma-1)}}+\frac{\gamma-1}{\gamma(3\gamma-1)}
z m_x)^2,\label{q eqn 3}
\end{eqnarray}
where we use $(m_x)^\backprime=-c m_{xx}$. After simplification, we
have \beq \tilde{q}^\backprime=\tilde{a}_0 -\tilde{a}_1 \tilde{q}+
\tilde{a}_2 \tilde{q}^2,\label{q fun}\eeq where $\tilde{a}_i$ are in
(\ref{tilde a0})$\sim$(\ref{tilde a2}). Similarly, we denote \beq
q={\bar{\mu}}\tilde{q},\label{def q}\eeq where ${\bar{\mu}}$ is in
(\ref{mu}). Since \beq {\bar{\mu}}^{\backprime}=\tilde{a}_1
{\bar{\mu}},\eeq (\ref{q fun}) changes to \beq q^\backprime=a_0+a_2
q^2,\eeq where $a_0$, $a_2$ are defined in (\ref{a0a2}).
\end{proof}

It is easy to check that $a_2<0$, \beq a_0\gtreqqless
0\Leftrightarrow (3\gamma-1)m m_{xx}\gtreqqless
(3\gamma+1)m_x^2\Leftrightarrow
(m^{-\frac{2}{3\gamma-1}})_{xx}\lesseqqgtr 0.\label{sign a0}\eeq So
the sign of $a_0$ only depends on the entropy profile, for fixed
$\gamma$. Note that the entropy profile of smooth solution is
stationary because of (\ref{lagrangian3 zm}). So the sign of $a_0$
only depends on the initial data.

{\cor{\label{cor2} When Assumption 1 holds,
\beq\left.\begin{array}{l} |y|\ \it{or}\ |q|\rightarrow \infty \
\it{iff}\ |u_x|\ \it{or}\ |\tau_x| \rightarrow
\infty.\end{array}\right.\label{y q blowup}\eeq}}\begin{proof} By
(\ref{def tilde y}), (\ref{def y}), (\ref{def tilde q}) and
(\ref{def q}),
\begin{eqnarray} y+q&=&2 m^{-\frac{3(3-\gamma)}{2(3\gamma-1)}}
u_x z^{\frac{\gamma+1}{2(\gamma-1)}},\\
y-q&=&2 m^{-\frac{3(3-\gamma)}{2(3\gamma-1)}}(m
z^{\frac{\gamma+1}{2(\gamma-1)}}
z_x+\frac{3\gamma-3}{3\gamma-1}z^{\frac{\gamma+1}{2(\gamma-1)}+1}m_x).\end{eqnarray}
By the same argument as the proof of (\ref{alpha beta blowup}) in
Lemma \ref{def2}, (\ref{y q blowup}) is right.  We complete the
proof.
\end{proof}
Similarly, when Assumption 1 holds, \beq\left.\begin{array}{l}
|\tilde{y}|\ \it{or}\ |\tilde{q}|\rightarrow  \infty \ \it{iff}\
|u_x|\ \it{or}\ |\tau_x| \rightarrow  \infty.\label{ty tq
blowup}\end{array}\right.\eeq
%
%
\section{Singularity Formation}
In this section, using the \emph{quadratic} ``ODEs" in the Theorem
\ref{them new ODEs}, we will give several singularity formation
(gradient blowup) results for
(\ref{lagrangian1})$\sim$(\ref{lagrangian3}) with smooth initial
data. By Corollary \ref{cor2}, we need to find conditions on the
initial data, under which $|y|$ or $|q|$ blows up in finite time.

\begin{figure}[htp] \centerline{
\includegraphics[scale=0.4]{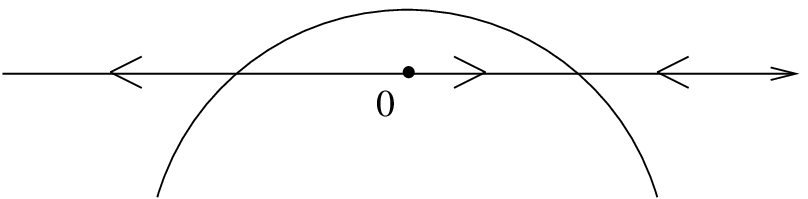}
\includegraphics[scale=0.4]{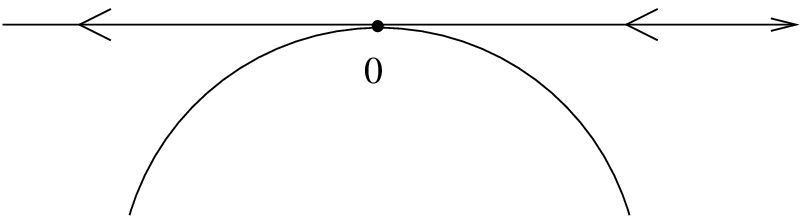}
\includegraphics[scale=0.4]{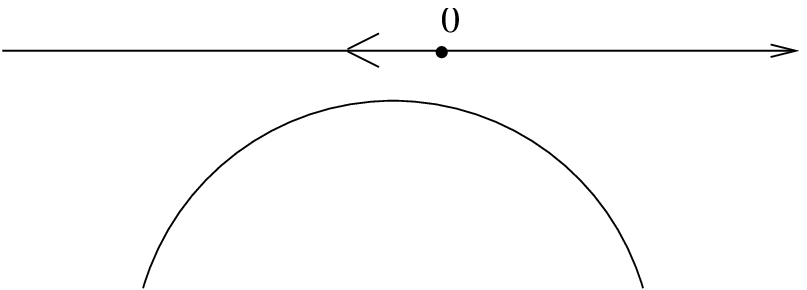}} \vspace*{8pt}
\caption{ $a_0>0$, $a_0=0$, $a_0<0$ from left to right. The arrows
indicate $y$ or $q$ increases or decreases.}
 \label{dynamic property}
\end{figure}

We consider the quadratic equation ($\xi= y$ or $q$), \beq 0=a_0+a_2
\xi^2, \label{quadricFun}\eeq where $a_0$ and $a_2$ are defined in
(\ref{a0def}) and (\ref{a2def}). When $a_0<0$, there is no real
root; when $a_0=0$, there is a unique root $0$; when $a_0>0$, the
two roots of (\ref{quadricFun}) are \beq\left.\begin{array}{l}
\pm\sqrt{
-\frac{a_0}{a_2}}=\pm\sqrt{\frac{2(\gamma-1)^2}{\gamma(\gamma+1)(3\gamma-1)}(m
m_{xx}-\frac{3\gamma+1}{3\gamma-1}m_x^2)}\
z^{\frac{\gamma+1}{2(\gamma-1)}+1}
m^{-\frac{3(3-\gamma)}{2(3\gamma-1)}}.\label{solutions}\end{array}\right.\eeq
By considering (\ref{new ode1}) and (\ref{new ode2}), we summarize
the dynamic system properties for smooth solutions in Fig.
\ref{dynamic property}. Moreover, $\{y>0,q>0\}$ is an invariant
domain for the region with initial entropy profile satisfying
$a_0\geqslant 0$, by Fig. \ref{dynamic property}.

For the $p$-system, $m$ is a constant. Hence, $a_0=0$, and \beq
y=s_x
m^{-\frac{3(3-\gamma)}{2(3\gamma-1)}}z^{\frac{\gamma+1}{2(\gamma-1)}},
\quad q=r_x
m^{-\frac{3(3-\gamma)}{2(3\gamma-1)}}z^{\frac{\gamma+1}{2(\gamma-1)}},\eeq
so $y$ and $s_x$ ($q$ and $r_x$) have same signs. Hence, the
$p$-system follows the middle picture of Fig. \ref{dynamic
property}. When $y$ or $q$ is negative somewhere, $y$ or $q$ will
approach negative infinity in finite time. When $y$ and $q$ are
nonnegative everywhere in the initial data, there is no such kind of
singularity, since $\{ y>0, q>0\}$ is an invariant domain. These
results are given by {Lax}, c.f. \cite{lax2}.

We don't expect such ``clean" results when the entropy profile is
varying, because the dynamic properties in this case are much more
complicated than those in a constant entropy domain, as we can see
from Fig. \ref{dynamic property}. In \cite{young blake 1} and a
forthcoming paper \cite{G1}, there are shock free examples with
compressive waves and entropy jumps. We also expect the smooth
examples including compressive waves.

Using the nonlinear \emph{Riccati} type equations (\ref{new ode1})
and (\ref{new ode2}), we can prove Theorems \ref{Thm singularity2}
and \ref{singularity} for breakdown of $y$ and $q$, which directly
generalize {Lax's} results in \cite{lax2}. This kind of formation of
singularity corresponds to the formation of a shock.

Theorem \ref{Thm singularity2} says that the shock free solutions
don't include ``strong" compressive waves, since $y$, $q$,
$\tilde{y}$, $\tilde{q}$ and $\alpha$, $\beta$ are larger than some
negative constants in shock free solutions. When $M_3$, $M_4$ are
close to zero, i.e. the variation of entropy profile is mild, $N$
and $\tilde{N}$ are close to zero, where $N$ and $\tilde{N}$ are
given at (\ref{Ndef}) and (\ref{tildeNdef}). The shock free
solutions should be ``almost rarefactive". Theorem \ref{Thm
singularity2} is quite sharp at this time, since it is consistent
with {Lax}'s singularity results in \cite{lax2}. In fact, when
entropy is a constant, $N=\tilde{N}=0$, Theorem \ref{Thm
singularity2} is the same as the singularity results for the
$p$-system in \cite{lax2}. The bounds $N$ and $\tilde{N}$ only
depend on $Z_U$, $M_1$, $M_2$, $M_3$, $M_4$. $N$ and $\tilde{N}$
don't depend on $Z_L$, so $Z_L$ can be arbitrarily small.
\vspace{.1cm}

In Theorem {\ref{singularity}}, we consider the singularity
formation for some special initial entropy profiles. Here, we give
an exact example with entropy profile satisfying condition
(\ref{condition1}) and bounded away from infinity and negative
infinity.
For example, \beq m(x)=\left\{\begin{array}{ll}f(x),&x\leqslant A,\\
(e^{-x}+1)^{-\frac{3\gamma-1}{2}},&x>A,\end{array}\right.\label{example}\eeq
satisfies (\ref{condition1}), where $f(x)$ is a smooth function that
makes $m(x)$ smooth and bounded away from zero and infinity. A large
class of entropy functions satisfies (\ref{condition1}).
\subsection{Proof of Theorem {\ref{Thm singularity2}}}
\begin{proof}
By (\ref{solutions}) and Assumption 2, the roots of
(\ref{quadricFun}), if these exist, have uniform lower bound $-N$,
where if $1<\gamma\leqslant 3$, \beq\left.\begin{array}{l}
N=(1+\varepsilon)\sqrt{\frac{2(\gamma-1)^2}{\gamma(\gamma+1)(3\gamma-1)}M_2
M_4}\ Z_U^{\frac{\gamma+1}{2(\gamma-1)}+1}
M_1^{-\frac{3(3-\gamma)}{2(3\gamma-1)}};\label{Ndef}
\end{array}\right.\eeq  if $\gamma> 3$, $M_1$ changes to $M_2$ in
(\ref{Ndef}). Here, $\varepsilon$ is a fixed positive constant,
which can be arbitrarily small. So \beq a_0+a_2
\frac{N^2}{(1+\varepsilon)^2}<0,\label{proof of main theorem}\eeq
since $a_2<0$. Recall, $m$ is stationary by (\ref{lagrangian3 zm}).

If $y<-N$ somewhere in the initial data, then by (\ref{proof of main
theorem}) and $a_2<0$, \beq a_0+a_2
\frac{y^2}{(1+\varepsilon)^2}<0.\eeq So (\ref{new ode1}) gives that
\beq y^{\prime}=a_0+a_2 y^2<(1-\frac{1}{(1+\varepsilon)^2}) a_2
y^2.\label{lastL1}\eeq Hence $y$ is decreasing at this time by
(\ref{a2def}), so $y$ is always less than $-N$ along forward
characteristic line. Furthermore, by (\ref{lastL1}), \beq
\frac{1}{y(t)}\geqslant{\frac{1}{y(0)}-\int_0^t
(1-\frac{1}{(1+\varepsilon)^2}) {a_2}dt},\label{SS9 1}\eeq where the
integral is along the characteristic. By (\ref{a2def}) and
Assumption 2, $a_2$ is negative and bounded below. So the right hand
side of (\ref{SS9 1}) approaches zero in finite time, which gives
that $y(t)$ approaches $-\infty$ in finite time. By (\ref{y q
blowup}), $|\tau_x|$ and/or $|u_x|$ blow up. \vspace{.1cm}

We have another version of the singularity formation results if we
start from (\ref{y fun}) and (\ref{q fun}) for $\tilde{y}$ and
$\tilde{q}$. We assume that\beq\left.\begin{array}{l}
\tilde{N}=(1+\varepsilon)\frac{(\gamma-1) [\
|9-3\gamma|M_3+\sqrt{|A_1|M_3^2+|A_2| M_2 M_4}\
]}{2(3\gamma-1)(\gamma+1)}{Z_U}^{\frac{\gamma+1}{2(\gamma-1)}+1},\label{tildeNdef}
\end{array}\right.\eeq where $\varepsilon$ is a fixed positive constant which
can be arbitrarily small, \beq\left.\begin{array}{lll} A_1&=&(9
\gamma^2-54\gamma+81)-\frac{1}{1+\varepsilon}\frac{1}{\gamma}(24\gamma^2+32\gamma+8),\\
A_2&=&\frac{1}{1+\varepsilon}\frac{1}{\gamma}(24\gamma^2+16\gamma-8).
\end{array}\right.\eeq The solutions of
\beq \tilde{a}_0\pm \tilde{a}_1
\tilde{y}+\frac{1}{1+\varepsilon}\tilde{a}_2 \tilde{y}^2=0,\eeq if
existing, are \beq\left.\begin{array}{l}
(1+\varepsilon)\Large{{\frac{({\gamma-1})
 [\ \pm{(9-3\gamma)} m_x\pm
\sqrt{A_1 m_x^2+ A_2 m m_{xx}}\
]z^{\frac{\gamma+1}{2(\gamma-1)}+1}}{2(3\gamma-1)(\gamma+1)}}},\label{proof
of main theorem 2}
\end{array}\right.\eeq
which are always larger than $-\tilde{N}$. If $\tilde{y}<-\tilde{N}$
somewhere in the initial data \beq
\tilde{y}^{\prime}<(1-\frac{1}{(1+\varepsilon)}) \tilde{a}_2
\tilde{y}^2.\eeq Hence $\tilde{y}$ will go to negative infinity in
finite time, by the same proof of the singularity result for $y$. By
(\ref{ty tq blowup}), the absolute values of $\tau_x$ and/or $u_x$
blow up.

We have symmetric results for $q$ and $\tilde{q}$.
\end{proof}
\subsection{Proof of Theorem \ref{singularity}}
\begin{proof}
We use $y(t)$ to denote the function $y$ along the forward
characteristic starting from $(0,A^*)$, and assume $y_0=y(0)<0$. By
(\ref{sign a0}), $a_0\leqslant0$ when $x>A$. So (\ref{new ode1})
changes to \beq y^{\prime}\geqslant a_2 y^2.\label{proof of second
thm}\eeq By solving (\ref{proof of second thm}), we get\beq
\frac{1}{y}\leqslant {\frac{1}{y_0}-\int a_2 dt},\label{proof of
second thm2}\eeq where the integral is along the characteristic.
Hence, when the right hand side of (\ref{proof of second thm2}) goes
to zero, $y$ goes to $-\infty$. So, by $y_0<0$, (\ref{a2def}) and
Assumption 2, the blowup happens before $T_*\leqslant-\frac{1}{y_0
\min(-a_2)}$, where $\min(-a_2)$ is a positive constant depending on
$Z_L$ ($Z_U$) and $M_1$ ($M_2$) when $1<\gamma<3$ ($\gamma>3$). When
$\gamma=3$, $a_2$ is a constant, so we don't need Assumption 2.
Symmetric results apply for $q$.
\end{proof}

When $1<\gamma<3$, we don't need the upper bounds for $z$, $|m_x|$
and $|m_{xx}|$.
\\{\bf Acknowledgement:} I am
grateful for the help of Professor {R.\ Young}, who leads me in this
area. He also gives me a lot of ideas and carefully revises this
paper. Professors {B.\ Temple} and {H.\ K.\ Jenssen} also gave me
some suggestions and inspirations.


\begin{thebibliography}{21}
\bibitem{Bressan} A. Bressan, \emph{Hyperbolic systems of conservation laws: the
one-dimmensional Cauchy problem}, Oxford Lecture Ser. Math. Appl.,
(Oxford Univ. Press, Oxford 2000).

\bibitem{G1} G. Chen and R. Young, Shock formation and exact solutions
for the compressible Euler equations with multiple entropy jumps,
\emph{in preparation}.

\bibitem{courant} R. Courant and K. O. Friedrichs, \emph{Supersonic flow
and shock waves}, (Wiley-Interscience, New York, 1948).

\bibitem{Dafermos} C. Dafermos, \emph{Hyperbolic conservation laws in
continuum physics}, Third edition, (Springer-Verlag, Heidelberg
2010).

\bibitem{Glimm} J. Glimm, Solutions in the large for nonlinear
hyperbolic systems of equations, \emph{Comm. Pure Appl. Math.}, {\bf
18} (1965) 697-715.

\bibitem{Glimm Lax} J. Glimm and P. Lax, Decay of solutions of
systems of nonlinear hyperbolic conservation laws, \emph{Amer. Math.
Soc. Memoir}, {\bf 101}. (Amer. Math. Soc: Providence, 1970).

\bibitem{Fritz John} F. John, Formation of singularities in one-dimensional
nonlinear wave propagation,  \emph{Comm. Pure Appl. Math.}, {\bf 27}
(1974) 377-405.

\bibitem{lax} P. Lax, Hyperbolic systems of conservation laws, II,
\emph{Comm. Pure Appl. Math.}, {\bf 10}, (1957) 537-566.

\bibitem{lax2} P. Lax, Development of singularities of solutions
of nonlinear hyperbolic partial differential equations, \emph{J.
Math. Physics}, {\bf 5:5}, (1964) 611-614.

\bibitem{lax3} P. Lax,  Hyperbolic systems of conservation laws and
the mathematical theory of shock waves, \emph{Conf. Board Math.
Sci.}, {\bf 11}, (SIAM, 1973).

\bibitem{Li daqian} T.-T. Li, Y. Zhou and D.-X. Kong, Global classical solutions
for general quasilinear hyperbolic systems with decay initial data,
\emph{Nonlinear Analysis, Theory, Methods $\&$ Applications}, {\bf
28:8}, (1997) 1299-1332.

\bibitem{Liu1} T.-P. Liu, The development of singularities in the
nonlinear waves for quasi-linear hyperbolic partial differential
equations, \emph{Jour. Diff. Equations}, {\bf 33}, (1979) 92-111.

\bibitem{Rammaha}  M. A. Rammaha, Formation of singularities in
compressible fluids in two-space dimensions, \emph{Proc. Amer. Math.
Soc.} {\bf 107:3}, (1989) 705-714.

\bibitem{sideris} T. Sideris, Formation of singularities
in three-dimensional compressible fluids, \emph{Commun. Math.
Phys.}, {\bf 101}, (1985) 475-485.

\bibitem{smoller} J. Smoller, \emph{Shock waves and reaction-diffusion
equations}, (Springer-Verlag, New York 1982).


\bibitem{young blake 1} B. Temple and R. Young,
A paradigm for time-periodic sound wave propagation in the
compressible Euler equations, \emph{Methods Appl. Anal.}, {\bf
16:3}, (2009) 341-364.

\bibitem{young blake 2} B. Temple and R. Young,
Time-periodic linearized solutions of the compressible Euler
equations and a problem of small divisors, \emph{to appear in SIAM
Jour. of Math. Anal.}.

\bibitem{wagner} D. Wagner,  Equivalence of the Euler and Lagrangian
equations of gas dynamics for weak solutions, \emph{Jour. Diff.
Equations}, {\bf 68}, (1987) 118-136.

\bibitem{young Global interaction} R. Young,
Global wave interactions in isentropic gas dynamics, \emph{to
appear}.

\bibitem{young com} R. Young, Convergence of characteristics and
shock formation, \emph{in preparation}.
\end{thebibliography}
\end{document}